\documentclass[10pt,twoside]{amsart}

\usepackage{amsthm}
\usepackage{amsmath}
\usepackage{amssymb}
\usepackage{amsfonts}
\usepackage{latexsym}
\usepackage{enumitem}
\usepackage{mathrsfs}
\usepackage[all]{xy} \SelectTips{eu}{}
\usepackage{hyperref}

\newtheorem*{intsetup*}{Setup}

\newtheorem*{intNotations*}{Notations and Conventions}

\newtheorem{intthm}{Theorem}[]

\newtheorem*{intque*}{Question}
\newtheorem*{intexa*}{Example}

\newcommand{\numberseries}{\bfseries}   
\newlength{\thmtopspace}                
\newlength{\thmbotspace}                
\newlength{\thmheadspace}               
\newlength{\thmindent}                  
\newtheoremstyle{bfupright head,slanted body}
                {\thmtopspace}{\thmbotspace}
                {\slshape}{\thmindent}{\bfseries}{.}{\thmheadspace}
                {{\numberseries \thmnumber{#2\;}}\thmnote{#3}}

\newtheoremstyle{bfupright head,upright body}
                {\thmtopspace}{\thmbotspace}
                {\upshape}{\thmindent}{\bfseries}{.}{\thmheadspace}
                {{\numberseries \thmnumber{#2\;}}\thmnote{#3}}

\newtheoremstyle{fixed bf head,slanted body}
                {\thmtopspace}{\thmbotspace}{\slshape}
                {\thmindent}{\bfseries}{.}{\thmheadspace}
                {{\numberseries \thmnumber{#2\;}}\thmname{#1}\thmnote{ (#3)}}

\newtheoremstyle{fixed bf head,upright body}
                {\thmtopspace}{\thmbotspace}{\upshape}
                {\thmindent}{\bfseries}{.}{\thmheadspace}
                {{\numberseries \thmnumber{#2\;}}\thmname{#1}\thmnote{ (#3)}}

\newtheoremstyle{numbered paragraph}
                {\thmtopspace}{\thmbotspace}{\upshape}
                {\thmindent}{\upshape}{}{\thmheadspace}
                {{\numberseries \thmnumber{#2.}}}

\theoremstyle{bfupright head,slanted body}
\newtheorem{res}{}[section]             \newtheorem*{res*}{}

\theoremstyle{bfupright head,upright body}
\newtheorem{bfhpg}[res]{}               \newtheorem*{bfhpg*}{}

\theoremstyle{fixed bf head,slanted body}
\newtheorem{thm}[res]{Theorem}          \newtheorem*{thm*}{Theorem}
\newtheorem{prp}[res]{Proposition}      \newtheorem*{prp*}{Proposition}
\newtheorem{cor}[res]{Corollary}        \newtheorem*{cor*}{Corollary}
\newtheorem{lem}[res]{Lemma}            \newtheorem*{lem*}{Lemma}
         \newtheorem*{que*}{Question}
\theoremstyle{fixed bf head,upright body}
\newtheorem{setup}[res]{Setup}           \newtheorem*{setup*}{Setup}
\newtheorem{dfn}[res]{Definition}       \newtheorem*{dfn*}{Definition}
\newtheorem{rmk}[res]{Remark}           \newtheorem*{rmk*}{Remark}
\newtheorem{exa}[res]{Example}           \newtheorem*{exa*}{Example}
           \newtheorem*{exer*}{Exercise}

\theoremstyle{numbered paragraph}
\newtheorem{ipg}[res]{}

\newlength{\thmlistleft}        
\newlength{\thmlistright}       
\newlength{\thmlistpartopsep}   
\newlength{\thmlisttopsep}      
\newlength{\thmlistparsep}      
\newlength{\thmlistitemsep}     

\newcounter{eqc}
  {\end{list}}%

\newcounter{prt}
\newenvironment{prt}{\begin{list}{\upshape (\alph{prt})}%
    {\usecounter{prt}%
      \setlength{\leftmargin}{\thmlistleft}%
      \setlength{\labelwidth}{\thmlistleft}%
      \setlength{\rightmargin}{\thmlistright}%
      \setlength{\partopsep}{\thmlistpartopsep}%
      \setlength{\topsep}{\thmlisttopsep}%
      \setlength{\parsep}{\thmlistparsep}%
      \setlength{\itemsep}{\thmlistitemsep}}}%
  {\end{list}}%

\newcounter{rqm}
\newenvironment{rqm}{\begin{list}{\upshape (\arabic{rqm})}%
    {\usecounter{rqm}%
      \setlength{\leftmargin}{\thmlistleft}%
      \setlength{\labelwidth}{\thmlistleft}%
      \setlength{\rightmargin}{\thmlistright}%
      \setlength{\partopsep}{\thmlistpartopsep}%
      \setlength{\topsep}{\thmlisttopsep}%
      \setlength{\parsep}{\thmlistparsep}%
      \setlength{\itemsep}{\thmlistitemsep}}}%
  {\end{list}}%

  {\end{list}}%

\newenvironment{prf*}[1][Proof]{%
  \begin{proof}[\bf #1]
    \setcounter{equation}{0}
    }
  {\end{proof}
}

\newcommand{\pgref}[1]{\ref{#1}}

\renewcommand{\eqref}[1]{(\pgref{eq:#1})}

\newcommand{\thmcite}[2][?]{\cite[Theorem~#1]{#2}}
\newcommand{\prpcite}[2][?]{\cite[Proposition~#1]{#2}}

\numberwithin{equation}{res}

\def\urltilda{\kern -.15em\lower .7ex\hbox{\~{}}\kern .04em}

\newcommand{\GP}{\mathsf{GP}}
\newcommand{\DP}{\mathsf{DP}}
\newcommand{\DI}{\mathsf{DI}}

\newcommand{\GI}{\mathsf{GI}}

\newcommand{\tr}{\mathsf{tr}}

\newcommand{\Proj}{\mathsf{Proj}}
\newcommand{\id}{\mathrm{id}}
\newcommand{\fd}{\mathrm{fd}}
\newcommand{\pd}{\mathrm{pd}}
\newcommand{\Inj}{\mathsf{Inj}}
\newcommand{\Flat}{\mathsf{Flat}}
\newcommand{\FI}{\mathsf{FI}}
\newcommand{\Mod}{\mathsf{Mod}}

\newcommand{\Coind}{\mathrm{Coind}}
\newcommand{\Ind}{\mathrm{Ind}}
\newcommand{\Id}{\mathrm{Id}}
\newcommand{\Rop}{R^{\sf op}}
\newcommand{\op}{\sf op}

\newcommand{\cok}{\mbox{\rm coker}}
\newcommand{\kernel}{\mbox{\rm ker}}
\newcommand{\im}{\mbox{\rm Im}}

\newcommand{\Hom}{\operatorname{Hom}}

\newcommand{\Ext}{\operatorname{Ext}}
\newcommand{\Tor}{\operatorname{Tor}}

\newcommand{\is}{\cong}

   \def\soft#1{\leavevmode\setbox0=\hbox{h}\dimen7=\ht0\advance
    \dimen7 by-1ex\relax\if t#1\relax\rlap{\raise.6\dimen7
    \hbox{\kern.3ex\char'47}}#1\relax\else\if T#1\relax
    \rlap{\raise.5\dimen7\hbox{\kern1.3ex\char'47}}#1\relax
    \else\if d#1\relax\rlap{\raise.5\dimen7\hbox{\kern.9ex
    \char'47}}#1\relax\else\if D#1\relax\rlap{\raise.5\dimen7
    \hbox{\kern1.4ex\char'47}}#1\relax\else\if l#1\relax
    \rlap{\raise.5\dimen7\hbox{\kern.4ex\char'47}}#1\relax
    \else\if L#1\relax\rlap{\raise.5\dimen7\hbox{\kern.7ex
    \char'47}}#1\relax\else\message{accent \string\soft
    \space #1 not defined!}#1\relax\fi\fi\fi\fi\fi\fi}

\begin{document}

\title[Relative Gorenstein modules over tensor rings]%
{$\mathcal{X}$-Gorenstein projective and $\mathcal{Y}$-Gorenstein injective modules over tensor rings}

\author[G.Q. Zhao]{Guoqiang Zhao}
\address{Guoqiang Zhao: Department of Mathmatics,
Hangzhou Dianzi University, Hangzhou 310018, China}
\email{gqzhao@hdu.edu.cn}

\author[J.X. Sun]{Juxiang Sun}
\address{Juxiang Sun (Corresponding author):
School of Mathematics and Statistics, Shangqiu Normal University, Shangqiu 476000, China}
\email{Sunjx8078@163.com}

\thanks{This research was partially supported by the National Natural Science Foundation of China (Grant No. 12061026) and the Key Research Projects Plan of Higher Education Institutions in Henan Province (Grant No. 26B110007).}





\begin{abstract}
Let $T_R(M)$ be a tensor ring and $\mathcal{X}$, $\mathcal{Y}$ be two classes of $R$-modules.
Under certain conditions, we prove that a $T_R(M)$-module $(A, u)$ is $\Ind(\mathcal{X})$-Gorenstein projective if and only if
$u$ is monomorphic and $\cok(u)$ is an $\mathcal{X}$-Gorenstein projective $R$-module.
$\mathcal{Y}$-Gorenstein injective $T_R(M)$-modules are also explicitly described.
As a consequence, the characterizations of Ding projective and Ding injective modules over $T_R(M)$ are obtained. 
Some applications to trivial ring extensions and Morita context rings are given.
\end{abstract}

\maketitle

\thispagestyle{empty}

\section*{Introduction}
\label{Preliminaries}

Throughout the paper, all rings are nonzero associative rings with identity and
all modules are unitary.
For a ring $R$, we adopt the convention that an $R$-module is a left $R$-module;
right $R$-modules are viewed as modules over the opposite ring $\Rop$.

The concept of the Gorenstein dimension of finitely generated modules over noetherian rings was introduced by Auslander and Bridger \cite{1969AB}.
Enochs and Jenda \cite{1995EnochsGP} extended their ideas and introduced the concepts of Gorenstein projective and Gorenstein injective modules over arbitrary rings.
Ding, Li and Mao considered two special cases of the Gorenstein projective and Gorenstein injective modules in \cite{DLM}, 
which were later called Ding projective and Ding injective modules by Gillespie \cite{G}.
It is shown that these modules 
are very useful to produce new model structures in the categories of modules.
In a series of papers \cite{A, GP, 2022XI, 2020GPTRI, 2022Mao, 2023Mao, TriExtGPMao} and so on,
these homological modules over triangular matrix rings,
trivial ring extensions and Morita context rings are also investigated.

Let $R$ be a ring and $M$ an $R$-bimodule.
Recall that the tensor ring $T_R(M)=\bigoplus_{i=0}^\infty M^{\otimes_Ri}$,
where $M^{\otimes_R0}=R$ and
$M^{\otimes_R(i+1)}=M\otimes_R(M^{\otimes_Ri})$ for $i \geqslant 0$.
Examples of tensor rings include but are not limited to
trivial ring extensions, Morita context rings, triangular matrix rings and so on.
The classical homological properties of the tensor ring $T_R(M)$ are studied
in \cite{1991Cohn, 2012Ample, 1975Roganov} and so on.

In \cite{GLS} Geiss, Leclerc and Schroer studied a certain kind of tensor ring, which is 1-Gorenstein, 
whose modules yield a characteristic-free categorification of the root system. 
Inspired by this work, Chen and Lu \cite{CL} recently studied the Gorenstein
homological properties of a tensor ring $T_R(M)$ for
an $N$-nilpotent $R$-bimodule $M$, where $R$ is a noetherian ring, and $M$ is finitely generated on both sides. 
They characterized Gorenstein projective $T_R(M)$-modules in terms of Gorenstein projective $R$-modules, 
and gave the relation between Gorenstein dimensions of rings $R$ and $T_R(M)$.
Later, Di, Liang et al. \cite{DL} extended and strengthened these results to the case that $R$ is an arbitrary associative ring 
and $M$ is not necessarily finitely generated $R$-bimodule, which recover many known results over trivial ring extensions, Morita context rings, triangular matrix rings.
A natural question arises:

\begin{que*} 
How to describe Ding projective and Ding injective modules over $T_R(M)$ for an $N$-nilpotent $R$-bimodule $M$?
\end{que*}

The main purpose of the present paper is to provide an answer for the above question and
further give some relevant applications.
In fact, we prove a more general result as follows, which recovers both the cases of Gorenstein modules and Ding modules; see Theorem 2.7 and Theorem 3.5.

\begin{intthm} \label{THM GP}

\begin{rqm}
\item Suppose that the $R$-bimodule $M$ is compatible and admissible with respect to $\mathcal{X}$.
Then 
$\Ind(\mathcal{X})$-$\GP(T_R(M))$ $= \Phi(\mathcal{X}$-$\GP(R))$.
\item Suppose that the $R$-bimodule $M$ is co-compatible and co-admissible with respect to $\mathcal{Y}$.
Then 
$\Coind(\mathcal{Y})$-$\GP(T_R(M)^{op})$ $= \Psi(\mathcal{Y}$-$\GP(R^{op}))$.
\end{rqm}
\end{intthm}

In the above theorem, for a subcategory $\sf{X}$ of $R$-modules
and a subcategory $\sf{Y}$ of $R^{\op}$-modules,
the symbols $\Phi(\sf{X})$ and $\Psi(\sf{Y})$
denote the subcategories of $\Mod(T_R(M))$ and $\Mod(T_R(M)^{\op})$, respectively,
which are defined as follows:
\begin{align*}
&\Phi({\sf{X}}) = \{(X,u) \in \Mod(T_R(M)) \,|\,
u \textrm{ is a monomorphism and } \cok(u) \in \sf{X}\}; \\
&\Psi({\sf{Y}}) = \{[Y,v] \in \Mod(T_R(M)^{\op}) \,|\,
v \textrm{ is an epimorphism and } \kernel(v) \in \sf{Y}\}.
\end{align*}

The following result gives an answer to the question above, see Corollary 2.12 and 3.10.

\begin{intthm} \label{}
\begin{rqm}
\item Suppose $\Tor_{\geqslant 1}^R(M, M^{\otimes_Ri}\otimes_R F)=0$
for each $F\in \Flat(R)$ and $i \geqslant 1$.
If $\fd_RM<\infty$ and $\fd_{\Rop}M<\infty$,
then $\DP(T_R(M)) = \Phi(\DP(R))$.
\item Let $T_R(M)$ be a right coherent ring and $M_{R}$ be finitely generated.
If $\Ext_{\Rop}^{\geqslant 1}(M$, $\Hom_{\Rop}(M^{\otimes_Ri}, E))=0$ for any $E\in\FI(R^{op})$ and $i\geqslant 1$,
$\fd_{R}M<\infty$and $\pd_{\Rop}M<\infty$,
then $\DI(T_R(M)^{\op}) = \Psi(\DI(\Rop))$.
\end{rqm}
\end{intthm}
An example is also given to show that there exists a bimodule $M$ over some algebra satisfying
all the conditions in the above theorem.

We finish the paper with two applications of the above theorems.
We not only reobtain the earlier results in this direction,
but also get some new conclusions.
We first study Ding projective (resp. injective)
modules over trivial extension of rings, 
then we study when a module over a Morita context ring is Ding projective (resp. injective).
The corresponding results present some new characterizations of Ding modules.

\section{Preliminaries}
\label{Preliminaries}
\noindent
In this section, we fix some notation, recall relevant notions
and collect some necessary facts.
Throughout the paper, we denote by $\Mod(R)$ the category of $R$-modules,
and by $\Proj(R)$ (resp., $\Inj(R)$, $\Flat(R)$) the subcategory of $\Mod(R)$
consisting of all projective (resp., injective, flat) $R$-modules. 
For an $R$-module $X$, denote by $\pd_{R}X$ (resp., $\id_{R}X$, $\fd_{R}X$)
the projective (resp., injective, flat) dimension of $X$, and
by $X^+$ the character module $\Hom_\mathbb{Z}(X, \mathbb{Q/Z})$.

\begin{bfhpg}[\bf Tensor rings]\label{Tensor rings}
Let $R$ be a ring and $M$ an $R$-bimodule.
We write $M^{\otimes_R 0} = R$ and $M^{\otimes_R (i+1)} = M \otimes_R (M^{\otimes_R i})$
for $i\geqslant 0$.
For an integer $N \geqslant 0$, recall that $M$ is said to be $N$-\emph{nilpotent} if
$M^{\otimes_R(N+1)}=0$.
For an $N$-nilpotent $R$-bimodule $M$, we denote by $T_R(M)=\bigoplus_{i=0}^N M^{\otimes_Ri}$
the \emph{tensor ring} with respect to $M$. It is easy to check that $T_R(M)^{\op}\is T_{\Rop}(M)$.

It follows from \cite{CL} that the category $\Mod(T_R(M))$ of $T_R(M)$-modules
is equivalent to the category $\Gamma$ whose objects are the pairs $(X, u)$,
where $X \in \Mod(R)$ and $u: M\otimes_R X \to X$ is an $R$-homomorphism,
and morphisms from $(X, u)$ to $(X', u')$ are
those $R$-homomorphisms $f \in \Hom_R(X, X')$ such that
$f \circ u = u' \circ (M\otimes f)$.
The category $\Mod(T_R(M)^{\op})$ of $T_R(M)^{\op}$-modules
is equivalent to the category $\Omega$ whose objects are the pairs $[Y, v]$,
where $Y \in \Mod(\Rop)$ and $v: Y \to \Hom_{\Rop}(M, Y)$ is an $\Rop$-homomorphism,
and morphisms from $[Y, v]$ to $[Y', v']$ are
those $\Rop$-homomorphisms $g \in \Hom_{\Rop}(Y, Y')$
such that $\Hom_{\Rop}(M, g) \circ v = v' \circ g$.

Unless otherwise specified, in the paper, we always view a  $T_R(M)$-module as a pair $(X,u)$
with $X \in \Mod(R)$ and $u \in \Hom_R(M \otimes_R X, X)$,
and view a $T_R(M)^{\op}$-module as a pair $[Y, v]$
with $Y \in \Mod(\Rop)$ and $v \in \Hom_{\Rop}(Y,\Hom_{\Rop}(M,Y))$. Note that a $T_R(M)^{\op}$-module can be equivalently viewed as a pair $(Y,\overline{v})$,
where $\overline{v} \in \Hom_{\Rop}(Y\otimes_RM, Y)$ is the adjoint morphism of $v$.

Let $(X, u) \overset{f} \longrightarrow (X', u') \overset{g} \longrightarrow (X'', u'')$ be a sequence in $\Mod(T_R(M))$.
We mention that $(X, u) \overset{f} \longrightarrow (X', u') \overset{g} \longrightarrow (X'', u'')$
is exact if and only if the underlying sequence
$X \overset{f} \longrightarrow X' \overset{g} \longrightarrow X''$
is exact in $\Mod(R)$.
\end{bfhpg}

\setup Throughout this paper, we always let $M$ denote an $N$-nilpotent $R$-bimodule.

\begin{bfhpg}[\bf The forgetful functor and its adjoint]
\label{The forgetful functor and its}
There exists a \emph{forgetful functor}
\begin{center}
$U: \Mod(T_R(M))\to \Mod(R)$,
\end{center}
which maps a  $T_R(M)$-module $(X,u)$ to the underlying $R$-module $X$.
Recall from \cite[Lemma 2.1]{CL} that $U$ admits a left adjoint
\begin{center}
$\Ind: \Mod(R) \to \Mod(T_R(M))$,
\end{center}
defined as follows:
\begin{itemize}
\item For an $R$-module $X$, define
$\Ind(X)=(\bigoplus_{i=0}^N(M^{\otimes_Ri}\otimes_RX),c_X)$,
where $c_X$ is an inclusion from
$M\otimes_R(\bigoplus_{i=0}^N(M^{\otimes_Ri}\otimes_RX))\cong
\bigoplus_{i=1}^N(M^{\otimes_Ri}\otimes_RX)$
to
$\bigoplus_{i=0}^N(M^{\otimes_Ri}\otimes_RX)$.
\item For an $R$-homomorphism $f: X \to Y$,
the $T_R(M)$-homomorphism $\Ind(f): \Ind(X) \to \Ind(Y)$ can be viewed as
a formal diagonal matrix with diagonal elements
$M^{\otimes_Ri}\otimes f$ for $0 \leqslant i\leqslant N$.
\end{itemize}
\end{bfhpg}

\begin{bfhpg}[\bf The stalk functor and its adjoint]
\label{The stalk functor and its adjoints}
There exists a \emph{stalk functor}
\begin{center}
$S: \Mod(R)\to \Mod(T_R(M))$,
\end{center}
which maps an $R$-module $X$ to the  $T_R(M)$-module $(X, 0)$.
The functor $S$ admits a left adjoint
$C: \Mod(T_R(M)) \to \Mod(R)$ defined as follows:
\begin{itemize}
\item For a  $T_R(M)$-module $(X,u)$, define $C((X,u))=\cok(u)$.
\item For a morphism $f: (X,u) \to (Y,v)$ in $\Mod(T_R(M))$, the morphism
$C(f) : \cok(u) \to \cok(v)$ is induced by the universal property of cokernels.
\end{itemize}

It is easy to see that $C\circ\Ind=\Id_{\Mod(R)}$, and
it follows from the Eilenberg-Watts theorem that the functor $\Ind$ is isomorphic to the functor $T_R(M)\otimes_R-$, so one has $\Ind(X)\is T_R(M)\otimes_RX$ for each $R$-module $X$.
\end{bfhpg}

\begin{ipg}\label{adjoint pairs-dual}
Similarly to \ref{The forgetful functor and its} and \ref{The stalk functor and its adjoints}, there are adjoint pairs $(S, K)$ and $(U, \Coind)$ as follows:
  \begin{equation*}
  \xymatrix@C=3pc{
    \Mod(\Rop)
    \ar[r]^-{S}
    &
    \Mod(T_R(M)^{\op})
    \ar@/^1.8pc/[l]_-{K}
    \ar[r]^-{U}
    &
    \Mod(\Rop).
    \ar@/^1.8pc/[l]_-{\Coind}
  }
  \end{equation*}
Here, $K([Y,v])=\kernel(v)$ for a $T_R(M)^{\op}$-module $[Y,v]$, and
$$\Coind(Y)=[\bigoplus_{i=0}^N\Hom_{\Rop}(M^{\otimes_Ri},Y),r_Y]$$
for an $\Rop$-module $Y$, where $r_Y$ is the morphism from
$\bigoplus_{i=0}^N\Hom_{\Rop}(M^{\otimes_Ri},Y)$ to
\[\Hom_{\Rop}(M,\bigoplus_{i=0}^N\Hom_{\Rop}(M^{\otimes_Ri},Y)) \cong
\bigoplus_{i=1}^N\Hom_{\Rop}(M^{\otimes_Ri},Y).\]
It is easy to see that $K\circ\Coind=\Id_{\Mod(\Rop)}$, and
it follows from the Eilenberg-Watts theorem that the functor $\Coind$ is isomorphic to the functor $\Hom_{\Rop}(T_R(M),-)$, so one has $\Coind(Y)\is \Hom_{\Rop}(T_R(M),Y)$ for each $\Rop$-module $Y$.
\end{ipg}

The following result is due to \cite[Lemma 1.9]{DL},
which will be used frequently in the sequel.

\begin{lem} \label{homological modules}
The following equalities hold.
\begin{prt}
\item $\Proj(T_R(M))=\Ind(\Proj(R))=\Phi(\Proj(R))$.
\item $\Inj(T_R(M)^{\op})=\Coind(\Inj(R^{\op}))=\Psi(\Inj(R^{\op}))$.
\item $\Flat(T_R(M))=\Phi(\Flat(R))$.
\end{prt}
\end{lem}

\section{$\mathcal{X}$-Gorenstein projective $T_R(M)$-modules}
\label{}
\noindent
We begin with the following definitions.

\begin{bfhpg}[\bf $\mathcal{X}$-Gorenstein projective modules]
\label{}
Let $\mathcal{X}$ be a class of $R$-modules that contains all projective $R$-modules.
Recall from \cite{BO} that an $R$-module $G$ is called {\it $\mathcal{X}$-Gorenstein projective}, if there exists
an exact sequence of projective $R$-modules 
$P^{\bullet} = \cdots\to P^{-1} \to P^{0}\to P^{1}\to\cdots$ such that $G\cong\im(P^{-1}\to P^{0})$ 
and $\Hom_{R}(P^{\bullet}, X)$ is exact whenever $X\in\mathcal{X}$.
The sequence $P^{\bullet}$ is called an {\it $\mathcal{X}$-complete projective resolution}.

We denote by $\mathcal{X}$-$\GP(R)$ the subcategory of $\mathcal{X}$-Gorenstein projective $R$-modules.
\end{bfhpg}

\begin{rmk}
\begin{prt}
\item[(1)] If $\mathcal{X}=\Proj(R)$, then $\mathcal{X}$-$\GP(R)$ is the category $\GP(R)$ of Gorenstein projective $R$-modules \cite{GHD}.
\item[(2)] If $\mathcal{X}=\Flat(R)$, then $\mathcal{X}$-$\GP(R)$ is the category $\DP(R)$ of Ding projective $R$-modules \cite{DLM, G}.
\end{prt}
\end{rmk}

Inspired by the definitions in \cite{DL}, we introduce the following two notions.

\begin{dfn} \label{dfn of compatible}
An $R$-bimodule $M$ is said to be \emph{compatible with respect to $\mathcal{X}$}
if the following two conditions hold
for an $\Ind(\mathcal{X})$-complete projective resolution $Q^\bullet$ of $T_R(M)$-modules,
where $\Ind(\mathcal{X})= \{ \Ind(X) | X\in \mathcal{X}\}$.
\begin{prt}
\item[(C1)]
The complex $M\otimes_RU(Q^\bullet)$ is exact.
\item[(C2)]
The complex $\Hom_{T_R(M)}(Q^\bullet,\Ind(M\otimes_RX))$ is exact
for each $X \in \mathcal{X}$.
\end{prt}
\end{dfn}

\begin{dfn}
An $R$-bimodule $M$ is said to be \emph{admissible with respect to $\mathcal{X}$}, 
if $$\Ext_R^1(G,M^{\otimes_Ri}\otimes_R X)=0=\Tor^R_1(M, M^{\otimes_Ri}\otimes_R G)$$
for each $G \in \mathcal{X}$-$\GP(R)$, $X \in \mathcal{X}$ and $i\geqslant 0$.
\end{dfn}

We mention that when $\mathcal{X}=\Proj(R)$, the above two definitions coincide with the ones in \cite{DL}.



In the following, we prove that if $M$ is compatible and admissible with respect to $\mathcal{X}$,
then one has $\Ind(\mathcal{X})$-$\GP(T_R(M)) = \Phi(\mathcal{X}$-$\GP(R))$.

\begin{lem} \label{GPTRM in Phi}
Suppose that the $R$-bimodule $M$ is compatible with respect to $\mathcal{X}$. Then there is an inclusion
$\Ind(\mathcal{X})$-$\GP(T_R(M)) \subseteq \Phi(\mathcal{X}$-$\GP(R))$.
\end{lem}

\begin{prf*}
Let $(A,u)$ be in $\Ind(\mathcal{X})$-$\GP(T_R(M))$ with
$$Q^\bullet:\ \xymatrix@C=0.5cm{
\cdots \to \Ind(P^{-1}) \ar[r]^-{d^{-1}} & \Ind(P^{0})
\ar[r]^-{d^0} & \Ind(P^{1}) \to \cdots}$$
an $\Ind(\mathcal{X})$-complete projective resolution such that $(A,u)\is\ker(d^0)$,
where each $P^i$ is in $\Proj(R)$.
Then there is an exact complex
$$\xymatrix@C=0.5cm{
\cdots \to
\bigoplus_{i=0}^N(M^{\otimes_Ri}\otimes_RP^{-1})
\ar[r]^{\,\,d^{-1}} & \bigoplus_{i=0}^N(M^{\otimes_Ri}\otimes_RP^0)
\overset{d^{0}}\to \cdots}$$
in $\Mod(R)$.
Since $M\otimes_RU(Q^\bullet)$ is exact by assumption, the complex
$$\xymatrix@C=1cm{\cdots \to \bigoplus_{i=1}^N(M^{\otimes_Ri}\otimes_RP^{-1})
\ar[r]^{\,\, M\otimes d^{-1}} &
\bigoplus_{i=1}^N(M^{\otimes_Ri}\otimes_RP^0) \overset{M\otimes d^{0}}\to \cdots}$$
is exact. Thus, one gets a commutative diagram with exact rows and columns:
$$
\xymatrix@C=0.7cm@R=0.7cm{
  0 \ar[r]^{}
  & \bigoplus_{i=1}^N(M^{\otimes_Ri}\otimes_RP^{-1}) \ar[d]_{M\otimes d^{-1}}
  \ar[r]^{c_{P^{-1}}}
  & \bigoplus_{i=0}^N(M^{\otimes_Ri}\otimes_RP^{-1})
    \ar[d]_{d^{-1}} \ar[r]^{\qquad \qquad \pi_{P^{-1}}}
  & P^{-1} \ar[d]_{} \ar[r]^{} & 0  \\
  0 \ar[r]^{}
  & \bigoplus_{i=1}^N(M^{\otimes_Ri}\otimes_RP^0) \ar[d]_{M\otimes d^0}
    \ar[r]^{c_{P^0}}
  & \bigoplus_{i=0}^N(M^{\otimes_Ri}\otimes_RP^0)
    \ar[d]_{d^0} \ar[r]^{\qquad \qquad \pi_{P^0}}
  & P^0 \ar[d]_{} \ar[r]^{} & 0  \\
  0 \ar[r]^{}
  & \bigoplus_{i=1}^N(M^{\otimes_Ri}\otimes_RP^1)
  \ar[r]^{c_{P^1}}
  & \bigoplus_{i=0}^N(M^{\otimes_Ri}\otimes_RP^1)
  \ar[r]^{\qquad \qquad \pi_{P^1}}
  & P^1 \ar[r]^{} & 0. }$$
This implies that $u$ is a monomorphism by the Snake Lemma.
Now we only need to prove that $\cok(u)\in\mathcal{X}$-$\GP(R)$.
It is sufficient to show that $\Hom_R(P^\bullet, X)$ is exact for each $X \in \mathcal{X}$,
where $P^\bullet$ is the third non-zero column in the above commutative diagram.
Notice that there exists a short exact sequence
$$\xymatrix@C=1cm{0 \to
\Ind(M\otimes_R X) \ar[r]^-{\phi_{(X,0)}} &
\Ind(X) \ar[r]^-{\eta_{(X,0)}} & (X,0) \to 0}$$
in $\Mod(T_R(M))$ by \cite[Remark 1.5]{DL}.
Applying the functor $\Hom_{T_R(M)}(Q^\bullet, -)$ to this sequence, one gets a short exact sequence
\[\xymatrix@C=1.5cm{ 0 \longrightarrow  (Q^\bullet,\Ind(M\otimes_R X))\ar[r]^-{(Q^\bullet,\phi_{(X,0)})} &
(Q^\bullet,\Ind(X))
 \ar[r]^-{(Q^\bullet,\eta_{(X,0)})} &(Q^\bullet,(X,0)) \longrightarrow 0}
\]
of complexes as each term of $Q^\bullet$ is projective.
Because $Q^\bullet$ is an $\Ind(\mathcal{X})$-complete projective resolution, $\Hom_{T_R(M)}(Q^\bullet,\Ind(X))$ is exact. By assumption, the complex $\Hom_{T_R(M)}(Q^\bullet, \Ind(M\otimes_RX))$ is also exact.
Thus, $\Hom_{T_R(M)}(Q^\bullet, S(X))=\Hom_{T_R(M)}(Q^\bullet, (X, 0))$ is exact, which yields that $\Hom_R(C(Q^\bullet), X)$ is exact since $(C, S)$ is an adjoint pair.
However, the complex $C(Q^\bullet)$ is nothing but $P^\bullet$
as $C\circ\Ind=\Id_{\Mod(R)}$.
Cosequently, $\Hom_R(P^\bullet, X)$ is exact, as desired.
\end{prf*}

\begin{lem} \label{Phi in GPTRM}
Suppose that the $R$-bimodule $M$ is admissible with respect to $\mathcal{X}$.
Then there is an inclusion $\Phi(\mathcal{X}$-$\GP(R)) \subseteq \Ind(\mathcal{X})$-$\GP(T_R(M))$.
\end{lem}

\begin{prf*}
Let $(A,u)$ be in $\Phi(\mathcal{X}$-$\GP(R))$.
Then there is an exact sequence in $\Mod(R)$
$$ 0 \to M\otimes_RA \overset{u}\longrightarrow A
\overset{\rho}\longrightarrow \cok(u) \to 0 \qquad(3.1)$$
with $\cok(u) \in \mathcal{X}$-$\GP(R)$.

\textbf{Step 1}:

Since $\cok(u) \in \mathcal{X}$-$\GP(R)$,
there exists an exact sequence
$$ 0\to\cok(u)\overset{\iota}\to P^{0} \overset{f^{0}}\to P^{1} \overset{f^{1}}\to \cdots \qquad(3.2)$$
in $\Mod(R)$, which is $\Hom_{R}(-, X)$-exact for any $X\in \mathcal{X}$.
Set $a_0= \iota \circ \rho : A \to P^{0}$, and hence,
we obtain an $R$-homomorphism
$M\otimes a_0= M\otimes_RA \to M\otimes_RP^{0}$.
Since $\Proj(R)\subseteq\mathcal{X}$,
then $\Ext_R^1(\ker(f^{i}),M^{\otimes_Rj}\otimes_RP)=0$ for every $P\in\Proj(R)$ and $j\geqslant 0$ by assumption.
Applying the functor $\Hom_R(-,M\otimes_RP^{0})$ to (3.1),
there exists $a_1\in \Hom_{R}(A, M\otimes_RP^{0})$ such that
$a_1 \circ u= M\otimes a_0$.
Applying the functor $\Hom_R(-,M\otimes_RM\otimes_RP^{0})$ to (3.1) again,
there exists $a_2\in \Hom_{R}(A, M\otimes_RM\otimes_RP^{0})$ such that
$a_2 \circ u= M\otimes a_1$.
We iterate the above argument to obtain
$$a_j: A \to M^{\otimes_Rj}\otimes_R P^{0}$$
such that $a_j \circ u= M\otimes a_{j-1}$ holds for $1 \leqslant j \leqslant N$.
Set $\widetilde{a}=(a_0, a_1, \cdots , a_N)^\tr$,
which is an $R$-homomorphism from
$A$ to $\bigoplus_{i=0}^N(M^{\otimes_Ri}\otimes_RP^{0})$.
It is routine to check that $\widetilde{a} \circ u = c_{P^{0}} \circ (M\otimes \widetilde{a})$.
This implies that $\widetilde{a}=(a_0, a_1, \cdots , a_N)^\tr$
forms a $T_R(M)$-homomorphism from $(A, u)$ to $\Ind(P^{0})$.
Thus one has the following exact commutative diagram
\begin{equation*}
\xymatrix@C=0.5cm@R=1cm{
  0 \ar[r] & M\otimes_RA \ar[d]_{M\otimes \widetilde{a}}\ar[r]^{u} & A \ar[d]_{\widetilde{a}}\ar[r]^{\rho} & \cok(u) \ar[d]_{\iota}\ar[r] & 0  \\
  0 \ar[r] & \bigoplus_{i=1}^N(M^{\otimes_Ri}\otimes_RP^{0}) \ar[r]^{c_{P^{0}}} & \bigoplus_{i=0}^N(M^{\otimes_Ri}\otimes_RP^{0}) \ar[r]^{\qquad \qquad (1,0,\cdots,0)} & P^{0} \ar[r] & 0}
\end{equation*}
for the sake of $(1,0,\cdots,0) \circ \widetilde{a}=a_0=\iota \circ \rho$, which yields that $C(\widetilde{a})=\iota$. 
Since $\Tor^R_1(M$, $M^{\otimes_Rj}\otimes_R\ker(f^{i}))=0$ for$j\geqslant 0$ by assumption,
hence $M^{\otimes_Rj}\otimes \iota$ is a monomorphism for each $0 \leqslant j \leqslant N$.
Note that $c_{P^{0}}$ is a split monomorphism,
one gets that $M^{\otimes_Rj}\otimes c_{P^{0}}$ is a monomorphism
for each $0 \leqslant j \leqslant N$.
Thus, it follows from \cite[Lemma 2.9]{DL} that $\widetilde{a}$ is a monomorphism, 
and hence there is a short exact sequence
$$ 0 \to (A, u) \overset{\widetilde{a}\,}\rightarrow \Ind(P^{0})\to (Z, v) \to 0$$
in $\Mod(T_R(M))$.

Next, we prove that $(Z, v) \in \Phi(\mathcal{X}$-$\GP(R))$.
Because the following diagram is commutative with exact rows and columns:
\begin{equation*} \label{factorization 1}
\xymatrix@C=1cm@R=0.5cm{
  & 0 \ar[d]^{} & 0 \ar[d]^{}  \\
  & M\otimes_R A \ar[d]_{u} \ar[r]^{M\otimes \widetilde{a} \,\,\, \qquad}
  & \bigoplus_{i=1}^N(M^{\otimes_Ri}\otimes_RP^{0}) \ar[d]_{c_{P^{0}}} \ar[r]^{}
  & M\otimes_R Z \ar[d]_{v} \ar[r]^{} & 0  \\
  0 \ar[r]^{}
  & A \ar[d]_{\rho} \ar[r]^-{\widetilde{a}}
  & \bigoplus_{i=0}^N(M^{\otimes_Ri}\otimes_RP^{0})
    \ar[d]_{\pi'} \ar[r]^{}
  & Z \ar[d]_{} \ar[r]^{} & 0  \\
  0 \ar[r]^{}
  & \cok(u) \ar[r]^{\quad \iota} \ar[d]^{}
  & P^{0} \ar[r]^{} \ar[d]^{}
  & \ker f^{1} \ar[r]^{}  & 0 , \\
  & 0  & 0   }
\end{equation*}
it induces an exact sequence
$0 \to M\otimes_R Z \overset{v \,}\rightarrow Z \to\ker f^{1} \to 0$
by the Snake Lemma, and so $(Z, v) \in \Phi(\mathcal{X}$-$\GP(R))$, as desired.

Continuing this process, we have the exact sequence of $T_R(M)$-modules
$$0 \to (A, u) \overset{\widetilde{a}}\to \Ind(P^{0})\overset{g^{0}}\to\Ind(P^{1}) \overset{g^{1}} \to \cdots. \qquad (3.3)$$

\textbf{Step 2}:

Since every $T_R(M)$-module is an image of a projective module, 
there exists a short exact sequence
$$ 0 \to (W, w) \overset{t}\longrightarrow \Ind(P^{-1})
             \overset{h}\longrightarrow (A, u) \to 0$$
of $T_R(M)$-modules with $P^{-1} \in \Proj(R)$.
Next, we show that $(W, w) \in \Phi(\mathcal{X}$-$\GP(R))$.
From the Snake Lemma,
there exists a short exact sequence
$$0 \to \cok(w) \to P^{-1} \overset{\pi}\to \cok(u) \to 0$$
of $R$-modules. Since $\cok(u)\in\mathcal{X}$- $\GP(R)$, then $\cok(w) \in \mathcal{X}$- $\GP(R)$ 
as $\mathcal{X}$- $\GP(R)$ is closed under taking kernels of epimorphisms by \cite[Theorem 2.3]{BO}.
Because $\Tor^R_1(M$, $M^{\otimes_Ri}\otimes_R\cok(u))=0$ for $0\leqslant i\leqslant N$ by assumption, 
it follows from \cite[Lemma 2.8]{DL} that $\Tor^R_1(M, A)=0$. Thus, $M\otimes t$ is a monomorphism, 
which yields that $w$ is a monomorphism since $t\circ w= c_{P^{-1}}\circ(M\otimes t)$.
Thus $(W, w) \in \Phi(\mathcal{X}$- $\GP(R))$.

Continuing this process, one gets the exact sequence in $\Mod(T_R(M))$
$$\cdots \to \rightarrow \Ind(P^{-2})\overset{g^{-2}}\to\Ind(P^{-1})  \overset{h}\to  (A, u)\to 0 \qquad (3.4)$$
and an exact sequence in $\Mod(R)$ $$\cdots\to P^{-2}\overset{f^{-2}} \to P^{-1} \overset{\pi}\to \cok(u) \to 0 \qquad (3.5)$$ 
Since $\cok(u) \in\mathcal{X}$- $ \GP(R)$, the sequence (3.5) is $\Hom_{R}(-, X)$-exact for any $X\in\mathcal{X}$.
Gluing (3.2) and (3.5) gives rise to an $\mathcal{X}$-complete projective resolution 
$$P^{\bullet}: \cdots\to P^{-2}\overset{f^{-2}}\to P^{-1}\overset{f^{-1}}\to P^{0} \overset{f^{0}}\to P^{1} \overset{f^{1}}\to \cdots$$
with $f^{-1}= \iota\circ\pi$ and $\ker f^{0} =\cok (u)$.

\textbf{Step 3}:

Combining (3.3) and (3.4) gives rise to an exact sequence of projective $T_R(M)$-modules
$$Q^{\bullet}: \cdots\to\Ind(P^{-2})\overset{g^{-2}}\to \Ind(P^{-1})\overset{g^{-1}}\to \Ind(P^{0})\overset{g^{0}}\to \Ind(P^{1})\overset{g^{1}}\to \cdots$$
with $g^{-1}= \widetilde{a}\circ h$ and $(A, u) = \ker g^{0}$. 
We claim that $Q^{\bullet}$ is $\Hom_{T_R(M)}(-, \Ind(X))$-exact for any $X\in\mathcal{X}$.

In fact, given $X\in\mathcal{X}$, 
since $(\Ind, U)$ is an adjoint pair, we have an isomorphism
$\Hom_{T_R(M)}(\Ind (P^{i}), \Ind(X))\cong$ 
$\Hom_{R}(P^{i}, \bigoplus_{i=0}^N(M^{\otimes_Ri}\otimes_R X))$. 
This implies that the complex 
\begin{center}
$\Hom_{T_R(M)}(Q^{\bullet}, \Ind(X))\cong$ 
$\bigoplus_{i=0}^N\Hom_{R}(P^{\bullet}, M^{\otimes_Ri}\otimes_R X)$
\end{center}
is exact since $\Ext_{R}^{1}(G, M^{\otimes_Ri}\otimes_R X) =0$ 
for any $G\in\mathcal{X}$-$\GP(R)$ and $i\geqslant 0$ by assumption.
\end{prf*}

Now the following result is obtained immediately from Lemma 2.5 and 2.6.

\begin{thm}\label{THM GP com}
Suppose that the $R$-bimodule $M$ is compatible and admissible with respect to $\mathcal{X}$.
Then there is an equality
$\Ind(\mathcal{X})$-$\GP(T_R(M))$ $= \Phi(\mathcal{X}$-$\GP(R))$.
\end{thm}

As a consequence, we can get the characterization of Gorenstein projective and Ding projective modules over $T_R(M)$.
The following lemma is needed.

\begin{lem} \label{homological modules}
 $\Flat(T_R(M))=\Ind(\Flat(R))=\Phi(\Flat(R))$.
\end{lem}

\begin{prf*}
The equality $\Flat(T_R(M))=\Phi(\Flat(R))$ holds by Lemma 1.6(c),
and it is clear that $\Ind(\Flat(R))\subseteq\Phi(\Flat(R))$. 
Thus, we have to prove the inclusion $\Flat(T_R(M)) \subseteq \Ind(\Flat(R))$.
Let $(A,u)\in\Flat(T_R(M))$. Then it is a direct limit of finitely generated free $T_R(M)$-modules
by \cite[Theorem 5.40]{Rot}. It follows from \cite[Lemma 1.7(c)]{DL} that $(A,u)\cong \underrightarrow{\lim}\Ind(F_{i})$, where each $F_{i}$ is finitely generated free $R$-module.
Since $(\Ind,U)$ is an adjoint pair, one has $\underrightarrow{\lim}\Ind(F_{i})\cong$ $\Ind(\underrightarrow{\lim}F_{i})$ from \cite[Theorem 5.43]{Rot}.
Note that $\underrightarrow{\lim}F_{i}$ is flat $R$-module by \cite[Theorem 5.40]{Rot} again.
Thus $(A,u)\in\Ind(\Flat(R))$.
\end{prf*}

The first part in the following result is the main result of \cite[Theorem 2.11]{DL}. However, the proof is a bit different.

\begin{cor}\label{G-proj}
\begin{prt}
\item
Suppose that the $R$-bimodule $M$  is compatible and admissible with respect to $\Proj(R)$. 
Then $\GP(T_R(M)) = \Phi(\GP(R))$.
\item
Suppose that the $R$-bimodule $M$  is compatible and admissible with respect to $\Flat(R)$. 
Then $\DP(T_R(M)) = \Phi(\DP(R))$.
\end{prt}
\end{cor}

\begin{prf*}
The assertion is obtained by Lemma 1.6(a), Theorem 2.7 and Lemma 2.8.
\end{prf*}

We give an example of an $R$-bimodule that is compatible and admissible with respect to $\Flat(R)$.
To this aim, we need the following lemma.

\begin{lem} 
If $\Tor_{\geqslant 1}^R(T_R(M), Z)=0$, then $\fd_R Z=\fd_{T_R(M)}\Ind(Z)$.
\end{lem}

\begin{prf*}
We may assume that $\fd_R Z= n$.
Then there is an exact sequence of $R$-modules
$0\to F_{n}\to \cdots\to F_{1}\to F_{0}\to Z\to 0$
with each $F_{i}$ flat. Applying $T_R(M)\otimes_{R}-$ to it gives rise to an exact sequence
$$0\to T_R(M)\otimes_{R}F_{n}\to \cdots\to T_R(M)\otimes_{R}F_{0}\to T_R(M)\otimes_{R}Z\to 0$$ 
by assumption. Note that $T_R(M)\otimes_{R}F_{i}=$ $\Ind F_{i}\in \Flat(T_R(M))$ by Lemma 2.8,
which means $\fd_{T_R(M)}\Ind(Z) \leqslant n$.

On the other hand, suppose $\fd_{T_R(M)}\Ind(Z)= m$. Take an exact sequence of $R$-modules
$0\to W\to F_{m-1}\cdots\to F_{1}\to F_{0}\to Z\to 0$ with each $F_{i}$ flat.
Applying $T_R(M)\otimes_{R}-$ to it gives rise to an exact sequence
$$0\to T_R(M)\otimes_{R}W\to T_R(M)\otimes_{R}F_{m-1}\to \cdots\to T_R(M)\otimes_{R}F_{0}\to T_R(M)\otimes_{R}Z\to 0.$$ 
Thus $T_R(M)\otimes_{R}W =\Ind(W)$ is flat as $\fd_{T_R(M)}\Ind(Z) = m$.
This infers $W\in\Flat(R)$ by Lemma 2.8 again, and hence $\fd_R Z \leqslant m$, we are done.
 \end{prf*}

\begin{prp} \label{example of admissible and compatible}
Suppose $\Tor_{\geqslant 1}^R(M, M^{\otimes_Ri}\otimes_R F)=0$
for each $F\in \Flat(R)$ and $i \geqslant 1$.
If $\fd_RM<\infty$ and $\fd_{\Rop}M<\infty$,
then $M$ is compatible and admissible with respect to $\Flat(R)$.
\end{prp}
\begin{prf*}
Let $Q^\bullet$ be a $\Flat(T_R(M))$-complete projective resolution.
Then $Q^\bullet$ can be written as:
$$\xymatrix@C=0.5cm{
\cdots \to \Ind(P^{-1}) \ar[r]^{\qquad d^{-1}\,\,} & \Ind(P^0)
\ar[r]^{d^0 \quad \,\,} & \Ind(P^1) \to \cdots,}$$
where $P^j \in \Proj(R)$ for each $j \in \mathbb{Z}$. Hence, the complex
$$U(Q^\bullet):\ \xymatrix@C=0.5cm{\cdots \to \bigoplus_{i=0}^N(M^{\otimes_Ri}\otimes_RP^{-1})
\ar[r]^{\,\, d^{-1}} &
\bigoplus_{i=0}^N(M^{\otimes_Ri}\otimes_RP^0) \to \cdots}$$ of $R$-modules is exact.
Next we prove that $M\otimes_R U(Q^\bullet)$ is exact.
Fix a flat resolution
$0 \to F^n \to \cdots \to  F^1 \to  F^0 \to M \to 0$ of $\Rop$-module $M$.
Because $\Tor_{\geqslant 1}^R(M,\bigoplus_{i=0}^N(M^{\otimes_Ri}\otimes_RP^j))=0$
by assumption, there exists an exact sequence of complexes
$$0 \to F^n\otimes_R U(Q^\bullet) \to \cdots \to  F^1\otimes_R U(Q^\bullet) \to
F^0\otimes_R U(Q^\bullet) \to M\otimes_R U(Q^\bullet) \to 0.$$
Since each $F^k\otimes_R U(Q^\bullet)$ is exact for $0\leqslant k\leqslant n$,
one gets that $M\otimes_R U(Q^\bullet)$ is exact, as desired.

We then prove the condition (C2). It follows from \cite[Lemma 2.4]{DL} that
$$\Tor_{\geqslant 1}^R(M^{\otimes_Ri},M\otimes_R F)=0$$
for each $F \in \Flat(R)$ and $1\leqslant i\leqslant N$,
and hence $\Tor_{\geqslant 1}^R(T_R(M),M\otimes_R F)=0$.
Therefore, by Lemma 2.10, one has $\fd_{T_R(M)}\Ind(M\otimes_R F)=\fd_R(M\otimes_R F)$.
Since $\fd_RM<\infty$, then $\fd_R(M\otimes_R F)<\infty$ by \cite[Lemma 3.2(2)]{2022Mao}. 
It follows from \cite[Lemma 3.1]{2022Mao} that $\Hom_{T_R(M)}(Q^\bullet,\Ind(M\otimes_RP))$ is exact.
Thus $M$ is compatible with respect to $\Flat(R)$.

Next, we prove that $M$ is admissible with respect to $\Flat(R)$.
Indeed, from \cite[Lemma 2.5]{DL} we know that $\fd_R(M^{\otimes_Rs})<\infty$ for each $0\leqslant s \leqslant N$ as $\fd_RM<\infty$. Thus, for each $F\in\Flat(R)$, $\fd_R(M^{\otimes_Rs}\otimes_R F)<\infty$ by \cite[Lemma 3.2(2)]{2022Mao}.
It follows from \cite[Lemma 3.1]{2022Mao} that $\Ext_R^{\geqslant 1}(G, M^{\otimes_Rs}\otimes_R F)=0$
for each $G \in \DP(R)$.
On the other hand, by another use of \cite[Lemma 2.5]{DL}, one has $\fd_{\Rop}(M^{\otimes_Rs})<\infty$,
which implies that $\Tor_{\geqslant 1}^R(M^{\otimes_Rs}, G)=0$ from \cite[Lemma 2.3]{EIT}.
Thus $\Tor_{\geqslant 1}^R(M,M^{\otimes_Rs}\otimes_RG)=0$ by \cite[Lemma 2.4]{DL}. 
\end{prf*}

\begin{cor} \label{example of admissible and compatible}
Suppose $\Tor_{\geqslant 1}^R(M, M^{\otimes_Ri}\otimes_R F)=0$
for each $F\in \Flat(R)$ and $i \geqslant 1$.
If $\fd_RM<\infty$ and $\fd_{\Rop}M<\infty$,
then $\DP(T_R(M)) = \Phi(\DP(R))$.
\end{cor}

\begin{prf*}
The conclusion follows from Corollary 2.9(b) and Proposition 2.11.
\end{prf*}

\section{$\mathcal{Y}$-Gorenstein injective $T_R(M)$-modules}
\label{}
\noindent

\begin{bfhpg}[\bf $\mathcal{Y}$-Gorenstein injective modules]
\label{}
Let $\mathcal{Y}$ be a class of right $R$-modules that contains all injective right $R$-modules.
Recall from \cite{GI} that
an exact complex $I^\bullet$ of injective right $R$-modules is called {\it $\mathcal{Y}$-complete injective resolution} if it remains exact after applying
the functor $\Hom_{R^{op}}(Y, -)$ for every $Y\in \mathcal{Y}$.
An right $R$-module $H$ is called \emph{$\mathcal{Y}$-Gorenstein injective}
provided that there exists a $\mathcal{Y}$-complete injective resolution $I^\bullet$
such that $H \cong \kernel(I^0 \to I^1)$. 
We denote by $\mathcal{Y}$-$\GI(R^{op})$ the subcategory of $\mathcal{Y}$-Gorenstein injective right $R$-modules.

Recall from \cite{S} that a right $R$-module $B$ is called $FP$-injective if $\Ext^{1}_{R}(N, B) = 0$ for every
finitely presented right $R$-module $N$, and
we denote by $\FI(R^{op})$ the subcategory of all $FP$-injective right $R$-modules. 
The $FP$-injective dimension of a right $R$-module $B$, denoted by $FP$-id$_{R^{op}}B$, is defined to be the smallest integer $n\geqslant 0$ such
that $\Ext^{n+1}_{R}(N, B) = 0$ for every finitely presented right $R$-module $N$ (if no such $n$
exists, set $FP$-$id_{R^{op}}B=\infty$).
\end{bfhpg}

\begin{rmk}
\begin{prt}
\item [(1)]If $\mathcal{Y}=\Inj(R^{op})$, then $\mathcal{Y}$-$\GI(R^{op})$ is the category $\GI(R^{op})$ of Gorenstein injective right $R$-modules \cite{GHD}.
\item[(2)] If $\mathcal{Y}=\FI(R^{op})$, then $\mathcal{Y}$-$\GI(R^{op})$ is the category $\DI(R^{op})$ of Ding injective right $R$-modules \cite{DLM, G}.
\end{prt}
\end{rmk}

\begin{dfn} \label{dfn of cocompatible}
An $R$-bimodule $M$ is said to be \emph{co}-\emph{compatible with respect to $\mathcal{Y}$} if
the following two conditions hold
for a $\mathcal{Y}$-complete injective resolution $I^\bullet$ of injective $T_R(M)^{\op}$-modules:
\begin{prt}
\item[(C1')] The complex $\Hom_{\Rop}(M, U(I^\bullet))$ is exact.
\item[(C2')] The complex $\Hom_{T_R(M)^{\op}}(\Coind(\Hom_{\Rop}(M, Y)), I^\bullet)$ is exact for each $Y\in\mathcal{Y}$.
\end{prt}
\end{dfn}

\begin{dfn} \label{dfn of coadmissible}
An $R$-bimodule $M$ is said to be \emph{co}-\emph{admissible with respect to $\mathcal{Y}$} if
\[
\Ext_{\Rop}^1(\Hom_{\Rop}(M^{\otimes_Ri}, Y), H)=0=\Ext_{\Rop}^1(M, \Hom_{\Rop}(M^{\otimes_Ri},H))
\]
for each $H \in$ $\mathcal{Y}$-$\GI(R^{op})$, $Y \in \mathcal{Y}$ and $i\geqslant 0$.
\end{dfn}


The next result is a dual version of Theorem 2.7, which can be proved dually.

\begin{thm} \label{THM GITRM = Psi cons}
Suppose that the $R$-bimodule $M$ is co-compatible and co-admissible with respect to $\mathcal{Y}$.
Then $\Coind(\mathcal{Y})$-$\GI(T_R(M)^{\op}) = \Psi(\mathcal{Y}$-$\GI(\Rop))$.
\end{thm}

To obtain the characterization of Ding injective modules over $T_R(M)$,
we need the following lemma.

\begin{lem} \label{homological modules}
Assume that $T_R(M)$ is a right coherent ring and $M_{R}$ is finitely generated. Then
\begin{prt}
\item[(1)] 
$\Coind(\FI(R^{op}))=\Psi(\FI(R^{op}))$.
\item[(2)] 
$\FI(T_R(M)^{op}) = \Psi(\FI(R^{op}))$.
\end{prt}
\end{lem}

\begin{prf*}
 Since $T_R(M)$ is a right coherent ring and $M_{R}$ is finitely generated, 
from \cite[Theorem 3.5]{DL} we know that $R$ is also right coherent and $M_{R}$ is finitely presented. Thus for any $B\in \Mod(R^{op})$, the natural map $\tau: M\otimes_{R} B^{+}\to \Hom_{R^{op}}(M, B)^{+}$ is an isomorphism by \cite[Lemma 3.55(ii)]{Rot}.

(1) It is cleat that $\Coind(\FI(R^{op}))\subseteq\Psi(\FI(R^{op}))$.
Given $[B, v]\in \Psi(\FI(R^{op}))$, there exits an exact sequence $0\to \ker v\to B\overset{v}\to\Hom_{R^{op}}(M, B)\to 0$
with $\ker v\in\FI(R^{op})$. This induces the exact sequence $$0\to\Hom_{R^{op}}(M, B)^{+}\overset{v^{+}}\to B^{+}\to (\ker v)^{+}\to 0,$$ where $(\ker v)^{+}\in\Flat(R)$ by \cite[Theorem 2.2]{Fi}.
Thus $(B^{+}, v^{+}\tau)\in\Phi(\Flat(R))=\Ind(F(R))$ by Lemma 2.8,
and hence $v$ is a split epimorphism with $\ker v\in \FI(R^{op})$.
By a dual argument of the proof of \cite[Lemma 2.3]{CL},
one gets $[B, v]\in\Coind(\FI(R^{op}))$, as desired.

(2) For $[B, v]\in\Mod(T_R(M)^{op})$, the sequence $0\to \ker v\to B\overset{v}\to\Hom_{R^{op}}(M, B)$ is exact 
if and only if the sequence $\Hom_{R^{op}}(M, B)^{+}\overset{v^{+}}\to B^{+}\to (\ker v)^{+}\to 0$ is exact.
Therefore, by \cite[Theorem 2.2]{Fi}, $[B, v]\in\FI(T_R(M)^{op})$ if and only if $[B, v]^{+}\cong(B^{+}, v^{+}\tau)$ $\in\Flat(T_R(M))$
$=\Phi(\Flat(R))$ from Lemma 1.6(c),
which means that $v^{+}$ is a monomorphism and $(\ker v)^{+}\cong$ $\cok v^{+}\cong\cok (v^{+}\tau)$$\in\Flat(R)$.
This is equivalent to that  $v$ is surjective and $\ker v\in\FI(R^{op})$  by \cite[Theorem 2.2]{Fi} again,
that is $[B, v]\in\Psi(\FI(R^{op}))$.
\end{prf*}

\begin{cor} \label{THM GITRM = Psi cons}
Ley $T_R(M)$ be a right coherent ring and $M_{R}$ be finitely generated.
If the $R$-bimodule $M$ is co-compatible and co-admissible with respect to $\FI(R^{op})$, 
then $\DI(T_R(M)^{\op}) = \Psi(\DI(\Rop))$.
\end{cor}

\begin{prf*}
The assertion follows from Theorem 3.5 and Lemma 3.6.
\end{prf*}

At the end of this section, we provide a situation that satisfies the conditions in Corollary 3.7.

\begin{lem} \label{Inj of Ind-}
Let $T_R(M)$ be a right coherent ring and $M_{R}$ be finitely generated.
Then for an $\Rop$-module $Y$, the following statements hold.
\begin{prt}
\item $Y$ is $FP$-injective if and only if $\Coind(Y)$ is an $FP$-injective $T_R(M)^{\op}$-module.
\item If $\Ext^{\geqslant 1}_{\Rop}(T_R(M),Y)=0$, then $FP$-$\id_{\Rop}Y=$ $FP$-$\id_{T_R(M)^{\op}}\Coind(Y)$.
\end{prt}
\end{lem}
 
 \begin{prf*}
It is an immediate consequence of Lemma 3.6.
\end{prf*}

\begin{prp} \label{example of coadmissible and cocompatible}
Let $T_R(M)$ be a right coherent ring and $M_{R}$ be finitely generated.
If $\Ext_{\Rop}^{\geqslant 1}(M, \Hom_{\Rop}(M^{\otimes_Ri}, E))=0$ for any $E\in\FI(R^{op})$ and $i\geqslant 1$,
$\fd_{R}M<\infty$and $\pd_{\Rop}M<\infty$,
then $M$ is co-compatible and co-admissible  with respect to $\FI(R^{op})$.
\end{prp}

\begin{prf*}
We first prove that $M$ is co-compatible with respect to $\FI(R^{op})$. Let
$$I^\bullet:\ \xymatrix@C=0.5cm{
\cdots \to \Coind(I^{-1}) \ar[r]^{\quad \,\,\,\,\, d^{-1}} & \Coind(I^0)
\ar[r]^{d^0 \quad \,\,} & \Coind(I^1) \to \cdots}$$
be a $\FI(T_R(M)^{op})$-complete injective resolution with each $I^{j}\in\Inj(R^{op})$. Then one gets an exact complex
$$U(I^\bullet):\ \xymatrix@C=0.5cm{ \cdots \to \bigoplus_{i=0}^N\Hom_{R^{\op}}(M^{\otimes_Ri},I^{-1})
\ar[r]^{\,\, d^{-1}} & \bigoplus_{i=0}^N\Hom_{R^{\op}}(M^{\otimes_Ri},I^{0}) \to \cdots}$$
of $\Rop$-modules.  We claim that the complex $\Hom_{\Rop}(M, U(I^\bullet))$ is exact. Fix a projective resolution $0 \to P^n \to \cdots \to  P^0 \to M \to 0$ of the $\Rop$-module $M$. By assumption one has $\Ext^{\geqslant 1}_{\Rop}(M,\bigoplus_{i=0}^N\Hom_{\Rop}(M^{\otimes_Ri}, I^{j}))=0$, so there is an exact sequence
$$0 \to \Hom_{\Rop}(M, U(I^\bullet)) \to \Hom_{\Rop}(P^0, U(I^\bullet))
\to \cdots \to \Hom_{\Rop}(P^n, U(I^\bullet)) \to 0$$
of complexes, which yields that the complex $\Hom_{\Rop}(M, U(I^\bullet))$ is exact, as claimed.

On the other hand, from \cite[Theorem 3.5]{DL} we know that $R$ is right coherent and $M_{R}$ is finitely presented. 
Thus for any $E\in\FI(R^{op})$, $\Ext_{\Rop}^{\geqslant 1}(M, E)=0$.
It follows from \cite[Lemma 2.15 (ii)]{DL} that
$$\Ext_{\Rop}^{\geqslant 1}(M^{\otimes_Rs},\Hom_{\Rop}(M,E))=0$$
for any $s\geqslant 1$, and so $\Ext_{\Rop}^{\geqslant 1}(T_R(M),\Hom_{\Rop}(M,E))=0$. 
Hence, Lemma 3.8 (b) gives rise to the equality
\begin{center}
$FP$-$\id_{T_R(M)^{\op}}\Coind(\Hom_{\Rop}(M,E))=$ $FP$-$\id_{\Rop}(\Hom_{\Rop}(M, E)),$
\end{center}
and it is finite by \cite[Lemma 4.2]{2022Mao} as $\fd_{R}M<\infty$.
It follows from \cite[Lemma 4.1]{2022Mao} that the complex 
$\Hom_{T_R(M)^{\op}}(\Coind(\Hom_{\Rop}(M,E)), I^\bullet)$ is exact, 
which yields that $M$ is co-compatible  with respect to $\FI(R^{op})$.

Next, we prove that $M$ is co-admissible  with respect to $\FI(R^{op})$.
Indeed, because for each $I \in \Inj(\Rop)$ and $i\geqslant 1$, 
$$\Hom_{\Rop}(\Tor_{\geqslant 1}^R(M, M^{\otimes_Ri}), I)\is\Ext^{\geqslant 1}_{\Rop}(M,\Hom_{\Rop}(M^{\otimes_Ri}, I))=0,$$
one has $\Tor_{\geqslant 1}^R(M, M^{\otimes_Ri})=0$. 
From \cite[Lemma 2.5]{DL} we know that $\fd_{R}(M^{\otimes_Ri})<\infty$ for each $0\leqslant i \leqslant N$, 
so $FP$-$\id_{\Rop}(\Hom_{\Rop}(M^{\otimes_Ri}, E))<\infty$
for each $E\in\FI(R^{op})$ by \cite[Lemma 4.2]{2022Mao} again,
which implies that $\Ext_{\Rop}^{\geqslant 1}(\Hom_{\Rop}(M^{\otimes_Ri}, E),H)=0$
for each $H \in \DI(\Rop)$. On the other hand,
by \cite[Lemma 2.5]{DL} again, one has $\pd_{\Rop}(M^{\otimes_Ri})<\infty$ for each $i\geqslant 0$,
which ensures that $\Hom_{\Rop}(M^{\otimes_Ri}, I^\bullet)$ is exact for an exact complex $I^\bullet$ 
consisting of injective $R^{op}$-modules by \cite[Lemma 2.5]{EIT}.
Therefore, $\Ext_{\Rop}^{\geqslant 1}(M^{\otimes_Ri}, H)=0$, which yields that
$\Ext^{\geqslant 1}_{\Rop}(M,\Hom_{\Rop}(M^{\otimes_Ri},H))=0$ for any $i\geqslant 0$
by \cite[Lemma 2.15 (i)]{DL}, as desired.
\end{prf*}

\begin{cor} 
Let $T_R(M)$ be a right coherent ring and $M_{R}$ be finitely generated.
If $\Ext_{\Rop}^{\geqslant 1}(M, \Hom_{\Rop}(M^{\otimes_Ri}, E))=0$ for any $E\in\FI(R^{op})$ and $i\geqslant 1$,
$\fd_{R}M<\infty$and $\pd_{\Rop}M<\infty$,
then $\DI(T_R(M)^{\op}) = \Psi(\DI(\Rop))$.
\end{cor} 

\begin{prf*}
The assertion follows from Corollary 3.7 and Proposition 3.9.
\end{prf*}

\section{Applications}
\label{Applications}
\noindent
In this section, we give some applications to
trivial ring extensions and Morita context rings.

\subsection{The trivial extension of rings}
\label{The trivial extension of rings}
\noindent
Let $R$ be a ring and $M$ an $R$-bimodule.
There exists a ring $R\ltimes M$,
where the addition is componentwise
and the multiplication is given by
$(r_1, m_1)(r_2, m_2) = (r_1r_2,r_1m_2 + m_1r_2)$ for $r_1, r_2 \in R$ and $m_1, m_2 \in M$.
This ring is called the \emph{trivial extension} of
the ring $R$ by the $R$-bimodule $M$; see \cite{TRIEXT1975} and \cite{TRIEXT1971}.

Suppose that the $R$-bimodule $M$ is $1$-nilpotent,
that is, $M\otimes_RM=0$.
Then it is easy to see that the tensor ring $T_R(M)$
is nothing but the trivial ring extension $R\ltimes M$.
One can immediately get the following results by
Theorem 2.7 and 3.5.

\begin{thm}\label{THM GP com}
Suppose that $M$ is a $1$-nilpotent $R$-bimodule.
\begin{rqm}
\item
If $M$ is compatible and admissible with respect to $\mathcal{X}$.
Then $\Ind(\mathcal{X})$-$\GP(R\ltimes M)$ $= \Phi(\mathcal{X}$-$\GP(R))$.
\item
Suppose that the $R$-bimodule $M$ is co-compatible and co-admissible with respect to $\mathcal{Y}$.
Then $\Coind(\mathcal{Y})$-$\GI(T_R(M)^{\op}) = \Psi(\mathcal{Y}$-$\GI(\Rop))$.
\end{rqm}
\end{thm}

The following result is obtained directly from Corollary 2.12 and 3.10; compare \cite[Theorem 2.1 and 3.2]{2023Mao}.

\begin{cor} \label{DPTRM = Phi in trivial extension}
Suppose that $M$ is a $1$-nilpotent $R$-bimodule.
\begin{rqm}
\item
If $\Tor_{\geqslant 1}^R(M,M\otimes_R F)=0$ for each $F\in \Flat(R)$,
$\fd_{R}M<\infty$ and $\fd_{\Rop}M<\infty$, then there is an equality $\DP(R\ltimes M) = \Phi(\DP(R))$. 

\item
Let $R\ltimes M$ be a right coherent ring and $M_{R}$ is finitely generated.
If $\Ext_{\Rop}^{\geqslant 1}(M$, $\Hom_{\Rop}(M, E))=0$ for any $E\in\FI(R^{op})$, 
$\fd_{R}M<\infty$ and $\pd_{\Rop}M<\infty$,
then there is an equality $\DI((R\ltimes M)^{\op}) = \Psi(\DI(\Rop))$.
\end{rqm}
\end{cor}

\subsection{Morita context rings}
\label{Morita context rings}
Let $A$ and $B$ be two rings,
and let $_AV_B$ and $_BU_A$ be two bimodules,
$\phi : U\otimes_AV \to B$ a homomorphism of $B$-bimodules,
and $\psi : V\otimes_BU \to A$ a homomorphism of $A$-bimodules.
Associated with a \emph{Morita context} $(A,B,U,V,\phi,\psi)$,
there exists a \emph{Morita context ring}
\[\Lambda_{(\phi,\psi)}=\begin{pmatrix}A & V \\ U & B\end{pmatrix};\]
see \cite{BASS1968, MORITA1958} for more details.
Following \thmcite[1.5]{GF1982}, one can view a $\Lambda_{(\phi,\psi)}$-module
as a quadruple $(X,Y,f,g)$ with $X \in \Mod(A)$, $Y \in \Mod(B)$,
$f \in \Hom_B(U\otimes_AX,Y)$, and $g \in \Hom_A(V\otimes_BY,X)$. Also, one can view a $\Lambda_{(\phi,\psi)}^{\op}$-module as a quadruple
$[W,N,s,t]$ with $W \in \Mod(A^{\op})$, $N \in \Mod(B^{\op})$,
$s \in \Hom_{B^{\op}}(N,\Hom_{A^{\op}}(U,W))$, and $t \in \Hom_{A^{\op}}(W,\Hom_{B^{\op}}(V,N))$.

It follows from \prpcite[2.5]{GFARTIN} that Morita context rings are trivial ring extensions
whenever both $\phi$ and $\psi$ are zero.
More precisely, consider the Morita context ring $\Lambda_{(0,0)}$.
There exists an isomorphism of rings:
$$\Lambda_{(0,0)}\overset{\cong}\longrightarrow (A\times B)\ltimes (U \oplus V)
\,\,\text{via}\,\, \begin{pmatrix}a & v \\ u & b\end{pmatrix} \mapsto
((a,b),(u,v)).$$
Thus, there exists an isomorphic functor
$$\mu: \Mod(\Lambda_{(0,0)}) \to \Mod((A\times B)\ltimes (U \oplus V))\ \mathrm{via}\
(X,Y,f,g) \mapsto ((X,Y),(g,f)),$$
where $(g,f)$ is from
$(U \oplus V)\otimes_{A\times B}(X, Y) \cong (V\otimes_BY,U\otimes_AX)$ to $(X, Y)$.

We mention that $(U \oplus V)\otimes_{A\times B}(U \oplus V) \cong (U\otimes_AV) \oplus (V\otimes_BU)$.
Then the $A\times B$-bimodule $U \oplus V$ is $1$-nilpotent if and only if $U\otimes_AV=0=V\otimes_BU$.
Thus, we obtain the next result by Corollary 4.2; compare \cite[Theorem 2.4 and 3.5]{2023Mao}.

\begin{cor} \label{con GP in morita ring}
Let $\Lambda_{(0,0)}$ be a Morita context ring with $U\otimes_AV=0=V\otimes_BU$, and let $(X,Y,f,g)$ be a $\Lambda_{(0,0)}$-module and $[W,N,s,t]$ a $\Lambda_{(0,0)}^{\op}$-module.
\begin{rqm}
\item
Suppose that $\Tor^B_{\geqslant 1}(V, U\otimes_A F_1)=0=\Tor^A_{\geqslant 1}(U, V\otimes_B F_2)$ for each $F_1 \in \Flat(A)$ and $F_2 \in \Flat(B)$. If $\fd_{B}U<\infty$, $\fd_{A^{\op}}U<\infty$, $\fd_{A}V<\infty$ and $\fd_{B^{\op}}V< \infty$, then $(X,Y,f,g)\in\DP(\Lambda_{(0,0)})$ if and only if
both $f$ and $g$ are monomorphisms, and $\cok(f)\in\DP(B)$ and $\cok(g)\in\DP(A)$. 
\item
Let $\Lambda_{(0,0)}$ be a right coherent ring and $U_{A}$, $V_{B}$ be finitely generated.
If $\Ext_{A^{op}}^{\geqslant 1}(U$, $\Hom_{B^{op}}(V, E_{2}))=0 $ 
$=\Ext_{B^{op}}^{\geqslant 1}(V$, $\Hom_{A^{op}}(U, E_{1}))=0 $
for any $E_{1}\in\FI(A^{op})$ and  $E_{2}\in\FI(B^{op})$, 
$\fd_{B}U<\infty$, $\fd_{A}V<\infty$, $\pd_{A^{op}}U<\infty$ and $\pd_{B^{op}}V<\infty$,
then $[W,N,s,t]\in\DI(\Lambda_{(0,0)}^{\op})$ if and only if both $s$ and $t$ are epimorphisms, and $\kernel(s)\in\DI(B^{\op})$ and $\kernel(t)\in\DI(A^{\op})$.
\end{rqm}
\end{cor}
\begin{prf*}
(1) Since $(U \oplus V)\otimes_{A\times B}(F_{1}, F_{2}) \cong (V\otimes_B F_{2},U\otimes_A F_{1})$,
then for each $(F_1, F_2) \in \Flat(A\times B)$, one has $\Tor^{A\times B}_{\geqslant 1}
(U \oplus V, (U \oplus V)\otimes_{A\times B}(F_1, F_2)) \cong$
$\Tor^{A}_{\geqslant 1}(U$, $V\otimes_B F_{2})\oplus$ $\Tor^{B}_{\geqslant 1}(V, U\otimes_A F_{1}) =0$
by assumption, and both $\fd_{A\times B}(U \oplus V)$ and $\fd_{{(A\times B)}^{\op}}(U \oplus V)$ are finite. 
Thus, it follows from Corollary 4.2(1) that there is an equality
$$
\DP(\Lambda_{(0,0)})
\cong \DP((A\times B)\ltimes (U \oplus V))
=\Phi(\DP(A\times B)),$$
where an object $((X,Y),(g,f))\in\Mod((A\times B)\ltimes (U \oplus V))$ is in $\Phi(\DP(A\times B))$ if and only if $(g,f)$ is a monomorphism and $\cok(g,f)$ is in $\DP(A\times B)$. This yields that $(X,Y,f,g)\in\DP(\Lambda_{(0,0)})$ if and only if
both $f$ and $g$ are monomorphisms, and $\cok(f)\in\DP(B)$ and $\cok(g)\in\DP(A)$.

(2) It follows from \cite[Corollary 4.7]{DL} that $A$ and $B$ are right coherent, 
and U and V are finitely presented as right $A$-module and $B$-module, respectively.
By a similar argument of (1), one can get the assertion.
\end{prf*}

\begin{rmk}
Corollary 4.3 improves the main result of \cite[Theorem 2.10]{A}, where the author only gave sufficient conditions for modules over $\Lambda_{(0,0)}$ to be Ding projective. 
\end{rmk}

We conclude this section with an example, due to \cite[Example 4.6]{DL}, showing that
there exists a bimodule $M$ over some algebra satisfying
the conditions in Corollary 4.3.

\begin{exa}\label{exa of motita UV=0}
Let $k$ be a field and $kQ$ the path algebra
associated to the quiver
\begin{align*}
\xymatrix{
Q:\  1 \ar[r]_{} & 2 \ar[r]_{} & 3 \ar[r]_{}
  & \cdots \ar[r]_{} & n \ar@/_3ex/[llll]_{}.}
\end{align*}
Suppose that $J$ is the ideal of $kQ$ generated by all the arrows.
Then $R=kQ/J^h$ is a self-injective algebra for $2\leqslant h \leqslant n$.
Denote by $e_i$ the idempotent element corresponding to the vertex $i$.
Then one has $e_jRe_i = 0$ whenever
$1\leqslant i< j\leqslant n$ and $j-i \geqslant h$.
Let $M = Re_i\otimes_k e_jR$. Then $M$ is an $R$-bimodule and projective on both sides, and $M\otimes_RM\cong Re_i\otimes_k (e_jR\otimes_RRe_i)\otimes_k e_jR=0$.
Thus, the Morita context ring
$$\Lambda_{(0,0)}=\begin{pmatrix}R & M \\ M & R\end{pmatrix}$$
satisfies all conditions in Corollary 4.3.
\end{exa}



\bibliographystyle{amsplain-nodash}

\begin{thebibliography}{10}

\bibitem{A}
D. Asefa, \emph{Ding projective modules over Morita context rings},
Comm. Algebra, {\bf 52} (2024), no. 1, 79--87. 

\bibitem{1969AB}
M. Auslander and M. Bridger, \emph{Stable module theory}, Memoirs of the
  American Mathematical Society, No. 94, American Mathematical Society,
  Providence, RI, 1969. 

\bibitem{BASS1968}
H. Bass, \emph{Algebraic {$K$}-theory}, W. A. Benjamin, Inc., New
  York-Amsterdam, 1968. 

\bibitem{BO}
D. Bennis and K. Ouarghi,  \emph{$\mathcal{X}$-Gorenstein projective modules}, International Mathematical Forum {\bf 5} (2010), no. 10, 487--491.

\bibitem{CL}
X.W. Chen and M. Lu, \emph{Gorenstein homological properties of tensor
  rings}, Nagoya Math. J. \textbf{237} (2020), 188--208. 

\bibitem{1991Cohn}
P.M. Cohn, \emph{Algebra. {V}ol. 3}, second ed., John Wiley \& Sons, Ltd.,
  Chichester, 1991. 

\bibitem{DL}
Z.X. Di, L. Liang, Z.Q. Song and G.L. Tang, \emph{Gorenstein homological modules over tensor rings}, arXiv:2504.21349.

\bibitem{DLM}
N.Q. Ding, Y.L. Li and L.X. Mao: Strongly Gorenstein flat modules. J. Aust. Math. Soc. 86 (2009), 323–338.

\bibitem{EIT}
E.E. Enochs, M.C. Izurdiaga and B. Torrecillas, \emph{Gorenstein conditions over triangular matrix rings}, J. Pure Appl. Algebra {\bf 218} (2014), 1544--1554.

\bibitem{1995EnochsGP}
E.E. Enochs and O.M.G. Jenda, \emph{Gorenstein injective and
  projective modules}, Math. Z. \textbf{220} (1995), no.~4, 611--633.


\bibitem{Fi}
D. J. Fieldhouse, \emph{Character modules, dimension and purity}, Glasg. Math. J. {\bf 13} (1972), 144--146.

\bibitem{TRIEXT1975}
R.M. Fossum, P.A. Griffith, and I. Reiten, \emph{Trivial
  extensions of abelian categories}, Lecture Notes in Mathematics, Vol. 456,
  Springer-Verlag, Berlin-New York, 1975, Homological algebra of trivial
  extensions of abelian categories with applications to ring theory.

\bibitem{GP}
N. Gao and C. Psaroudakis, \emph{Gorenstein homological aspects of
  monomorphism categories via {M}orita rings}, Algebr. Represent. Theory
  \textbf{20} (2017), no.~2, 487--529. 

\bibitem{GLS}
C. Geiss, B. Leclerc and J. Schroer, \emph{Quivers with relations for symmetrizable Cartan
matrices I: Foundations}, Invent. Math. \textbf{209} (2017), 61--158.

\bibitem{G}
J.Gillespie, \emph{Model structures on modules over Ding-Chen rings}, Homology Homotopy Appl. {\bf 12} (2010), 61--73.

\bibitem{GI}
J.Gillespie and A. Icaob, \emph{Duality pairs, generalized Gorenstein modules, and Ding injective envelopes},
Comptes Rendus. Mathématique {\bf 360} (2022), 381--398.

\bibitem{GF1982}
E.L. Green, \emph{On the representation theory of rings in matrix form},
  Pacific J. Math. \textbf{100} (1982), no.~1, 123--138. 

\bibitem{GFARTIN}
E.L. Green and C. Psaroudakis, \emph{On {A}rtin algebras arising
  from {M}orita contexts}, Algebr. Represent. Theory \textbf{17} (2014), no.~5,
  1485--1525. 

\bibitem{2022XI}
Q.Q. Guo and C.C. Xi, \emph{Gorenstein projective modules over rings
  of morita contexts}, Sci. China Math. \textbf{67} (2024), no.~11, 2453--2484. 

\bibitem{GHD}
H. Holm, \emph{Gorenstein homological dimensions}, J. Pure Appl. Algebra
  \textbf{189} (2004), no.~1-3, 167--193. 



\bibitem{2020GPTRI}
H.H. Li, Y.F. Zheng, J.S. Hu, and H.Y. Zhu, \emph{Gorenstein
  projective modules and recollements over triangular matrix rings}, Comm.
  Algebra \textbf{48} (2020), no.~11, 4932--4947. 


\bibitem{2022Mao}
L.X. Mao, \emph{Ding modules and dimensions over formal triangular matrix rings}, Rend. Semin. Mat. Univ. Padova {\bf 148} (2022), 1--22.

  
\bibitem{2023Mao}
\bysame, \emph{Ding projective and Ding injective modules over trivial ring extensions}, Czechoslovak Math. J. {\bf 73} (2023), no. 148, 903--919.  

\bibitem{TriExtGPMao}
\bysame, \emph{Gorenstein projective, injective and flat modules over trivial
  ring extensions}, J. Algebra Appl. \textbf{24} (2025), no.~2, 2550030, 23. 

\bibitem{2012Ample}
H. Minamoto, \emph{Ampleness of two-sided tilting complexes}, Int. Math.
  Res. Not. IMRN (2012), no.~1, 67--101. 

\bibitem{MORITA1958}
K. Morita, \emph{Duality for modules and its applications to the theory of
  rings with minimum condition}, Sci. Rep. Tokyo Kyoiku Daigaku Sect. A
  \textbf{6} (1958), 83--142. 

\bibitem{TRIEXT1971}
I. Reiten, \emph{Trivial extensions and Gorenstein rings}, ProQuest
  LLC, Ann Arbor, MI, 1971, Thesis (Ph.D.)--University of Illinois at
  Urbana-Champaign. 

\bibitem{1975Roganov}
J.V. Roganov, \emph{The dimension of the tensor algebra of a projective
  bimodule}, Mat. Zametki \textbf{18} (1975), no.~6, 895--902. 

\bibitem{Rot}
J. J. Rotman,  \emph{An Introduction to Homological Algebra}, 2nd ed., Springer Press, New York, 2008.

\bibitem{S}
B. Stenstr\"{o}m, \emph{Coherent rings and FP-injective modules}, J. London Math. Soc. {\bf 2} (1970),
323--329.




\end{thebibliography}

\def\cprime{$'$}
  \providecommand{\arxiv}[2][AC]{\mbox{\href{http://arxiv.org/abs/#2}{\sf
  arXiv:#2 [math.#1]}}}
  \providecommand{\oldarxiv}[2][AC]{\mbox{\href{http://arxiv.org/abs/math/#2}{\sf
  arXiv:math/#2
  [math.#1]}}}\providecommand{\MR}[1]{\mbox{\href{http://www.ams.org/mathscinet-getitem?mr=#1}{#1}}}
  \renewcommand{\MR}[1]{\mbox{\href{http://www.ams.org/mathscinet-getitem?mr=#1}{#1}}}
\providecommand{\bysame}{\leavevmode\hbox to3em{\hrulefill}\thinspace}
\providecommand{\MR}{\relax\ifhmode\unskip\space\fi MR }
\providecommand{\MRhref}[2]{%
  \href{http://www.ams.org/mathscinet-getitem?mr=#1}{#2}
}
\providecommand{\href}[2]{#2}

\end{document}